\documentclass[]{article}
\usepackage{amsmath,amsfonts,amssymb}
\usepackage{color}
\usepackage{graphicx}

\newtheorem{thm}{Theorem}[section]
\newtheorem{rem}{Remark}[section]
\title{Numerical solution to a Parabolic-ODE Solow model with spatial diffusion and technology-induced motility}
\author{N. Ure\~{n}a, A.M. Vargas}

\begin{document}

\maketitle

\begin{abstract}
This work studies a parabolic-ODE PDE's system which describes the evolution of the physical capital ``$k$"  and technological progress ``$A$",
using a meshless in one and two dimensional bounded domain  with regular boundary. The well-known Solow model is extended by considering the spatial diffusion of both capital anf technology. Moreover, we study the case in which no spatial diffusion of the technology progress occurs. For such models, we propound schemes based on the Generalized Finite Difference method and proof the convergence of the numerical solution to the continuous one. Several examples show the dynamics of the model for a wide range of parameters. These examples illustrate the accuary of the numerical method.
\end{abstract}

	\section{Solow model: Mathematical formulation}
The Solow model is an important theoretical framework in economics because it provides insights into the long-run growth of an economy. The model was developed by Robert Solow in the 1950s and 1960s \cite{solow}, and it is named after him. It is an important tool for economists to understand the determinants of long-term economic growth and to analyze the impact of various government policies on economic performance. The model explains how changes in population, capital accumulation, and technological progress affect economic growth. In particular, it highlights the role of technological progress in driving long-term economic growth.
\\
However, it has become clear that incorporating spatial considerations is crucial for a more comprehensive analysis.

The spatial dimension introduces an additional layer of complexity to economic growth dynamics. Economic activities, including the diffusion of knowledge, investment, and productivity, are not uniformly distributed across regions. Certain areas may benefit from agglomeration effects, technological spillovers, or specialized industries, leading to higher levels of productivity and growth. At the same time, regions with limited access to resources or infrastructure may face economic challenges and slower development.
By incorporating spatial diffusion into the Solow model, we can capture these spatial dynamics and gain deeper insights into the patterns of economic growth. Spatial diffusion allows us to examine how technological advancements, capital flows, and knowledge spillovers spread across regions, influencing the productivity levels and growth rates of different areas. This perspective enables us to understand the interplay between regional disparities, agglomeration effects, and the overall economic performance of a country or region. Models including spatial diffusion are \cite{gonzales}, \cite{grasseti}, \cite{juchem} and \cite{juchem2}.

Moreover, considering spatial diffusion in the Solow model can help policymakers design more effective regional development strategies. By understanding the mechanisms behind spatial disparities, governments and policymakers can identify areas that require targeted investments, infrastructure improvements, or policies to enhance knowledge diffusion and promote economic convergence
In this paper we model the second case, that is to say, we are interested in the situation where technology progress is constant or created at a positive rate but because of monopolistic or autarkic scenarios. 
Finally, our model reads as follows
\begin{equation}
	\label{1}
	\left\lbrace
	\begin{aligned}{}
		&\frac{\partial{k}}{\partial{t}} = \Delta{k} - div(\chi k \nabla{A}) + Af(k)-\delta k ,& x\in \Omega, \quad  t>0, \\
		&\frac{\partial{A}}{\partial{t}}=Ag_A, & x\in \Omega, \quad t>0, \\
		&k(x,0)=k_0(x), \qquad  A(x,0)=A_0(x), &  x\in \Omega,  \\
		&\frac{\partial{k}}{\partial{\nu}}=\frac{\partial{A}}{\partial{\nu}}=0, &  x\in \partial \Omega, \quad t>0,
	\end{aligned}
	\right.
\end{equation}
Here, $f$ denotes some production function and $\delta$ represents the depreciation rate (assumed constant). Several production functions are found in the literature, for instance the Cobb-Douglas prodcution function, $f(k)=k^\alpha.$ Due to the reasons exposed in \cite{capasso}, we use $$f(k)=\frac{\alpha_1 k^{p}}{1+\alpha_2 k^q}.$$ In \cite{capasso} and \cite{boqueron}, the authors assumed that the capital flows from regions with abundant capital toward the ones with relatively less capital. We can assume that the mobility of the capital is also induced by the technological progress and not only by random motion.
In this model, regions with higher levels of technology exhibit higher productivity and attract capital flows from regions with lower technology concentrations. This reflects the idea that capital tends to move towards areas where technological advancements offer greater economic opportunities. The model considers a non-concave production function to capture more realistic economic conditions, allowing for diminishing returns to capital accumulation. By incorporating capital flows towards technology concentration, the model explores how the spatial distribution of technology and capital affects long-term economic growth patterns.
\\
Numerical simulations offer a quantitative exploration of the model's dynamics, illustrating the effects of technology concentration and capital flows on economic growth trajectories. 	Understanding the implications of spatial heterogeneity in technology concentration and capital movements can provide valuable insights for policymakers and researchers interested in regional development strategies. By examining the interplay between technological advancements and capital flows, this study contributes to the broader understanding of spatial economics and the determinants of regional growth.
\\
To obtain empirical evidence of the dynamics within this nonlinear PDE system, we conducted long-term simulations using the Generalized Finite Difference Method (GFDM), a meshless approach based on truncated Taylor series and moving least squares. The GFDM has garnered significant attention since the influential work of Lizska and Orkisz \cite{li} and subsequent studies by Benito, Gavete, and Ure\~na \cite{R1}. The mesh-independent nature of the GFDM allows for the computation of numerical solutions in highly complex domains with irregular node distributions. This versatility renders the method a powerful tool for tackling real-world problems involving nonlinear PDEs.

The applications of GFDM span diverse domains, showcasing its efficacy. For instance, it has been successfully employed in studying chemotaxis systems in biology \cite{bbg}, elastic wave propagation problems \cite{salete}, and porous flow and geomechanics \cite{liu}. Detailed accounts of these applications can be found in the corresponding references, offering a comprehensive overview.

The consistency of the GFDM's explicit formulae has been established for one-dimensional cases \cite{vargas} as well as two-dimensional cases \cite{R6}. Therefore, for brevity, we omit a discussion of this aspect here. The primary focus of this paper is on Theorem \ref{thm}, where we demonstrate the conditional convergence of the GFDM's numerical solution. Specifically, as time (t) increases, the error between the numerical solution and the continuous solution of the model approaches zero, given that certain conditions pertaining to the time increment and scheme coefficients are satisfied.
\\
The  paper is organized as follows: in Section 2 we introduce  some  explicit formulae using the Generalized Finite-Difference method for 1d and 2d problems. Next, in Section 3 we study the convergence of the GFD explicit scheme and we prove Theorem \ref{thm}. In Section 4, extensive numerical experiments (convergence studies, long-time simulations, etc.) are presented to illustrate the accuracy, efficiency and robustness of the developed numerical algorithms. We finally present some conclusions.
\section{1D explicit formulae}
Consider a discretization $M$ of $[0,L]$ with $N$ nodes and a subset of $M$ of $s$ points $$\{x_i:i=1,...,s\} \subset M$$ with center at some node, $x_0$, called star. Different criteria to choose the $s$ nodes of the star can be found in \cite{R1} and \cite{lancaster}.
To find the discretization of spatial derivatives at each of the points of that star of some regular enough function $u$, we take the truncated second order Taylor expansion of the solution around $x_0$
\begin{equation}
	u(x_i)=u(x_0)+ \frac{\partial u}{\partial x}(x_0)\cdot(x_i-x_0)+\frac{1}{2}\frac{\partial^2 u}{\partial x^2}(x_0)\cdot(x_i-x_0)^2+\mathcal{O}(3).
\end{equation}
We denote as $U_i=U(x_i)$ the approximation of the $u$ at $x_1$. Further, we define  the vectors $$\textbf{d}=\Biggl(\frac{\partial U_i}{\partial x},\frac{\partial^2 U_i}{\partial x^2}\Biggr),\quad\textbf{c}_i=\displaystyle\Biggl(h_i,\frac{h_i^2}{2}\Biggr),$$ where  $h_i=x_i-x_0$, and the operator $B$ be the sum of the weighted quadratic errors (properties of the weighting functions can be found in \cite{lancaster}, Section 10.3 Moving least squares methods)
$$B(\textbf{d})=\sum_{i=1}^{s}(U_0-U_i+\textbf{c}_i^{T}\textbf{d})^2w_i^2.$$
By minimizing $B$ with respect to $\textbf{d}$ we arrive to the linear system
$$\sum_{i=1}^{s}w_i^2\textbf{c}_i\textbf{c}_i^{T}\textbf{d}=-\sum_{i=1}^{s}w_i^2(U_0-U_i)\textbf{c}_i.$$
As proved in \cite{R6}, the matrix $\Lambda:=\displaystyle\sum_{i=1}^sw_i\textbf{c}_i\textbf{c}_i^{T}$ is positive definite, thus
$$\textbf{d}=-U_0\sum_{i=1}^{s}w_i^2\Lambda^{-1}\textbf{c}_i+\sum_{i=1}^{s}U_iw_i^2\Lambda^{-1}\textbf{c}_i.$$
For simplicity, we define the vectors
$$\boldsymbol\lambda_0:=\sum_{i=1}^{s}w_i^2\Lambda^{-1}\textbf{c}_i,\quad \boldsymbol\lambda_i:=w_i^2\Lambda^{-1}\textbf{c}_i.$$
Then, introducing the notation $\boldsymbol\lambda_0=(\lambda_{01},\lambda_{02})^{T}$ (analogously for $\boldsymbol\lambda_i$) and the condition $\boldsymbol\lambda_0=\displaystyle\sum_{i=1}^s\boldsymbol\lambda_i$, we write the spatial derivatives of the function, by means of the components of the vectors $\boldsymbol\lambda_0$, $\boldsymbol\lambda_i$, as a linear combination of the values of the solution at the surrounding nodes:
\begin{equation}\label{discretizacion}
	\left\lbrace
	\begin{aligned}
		&\frac{\partial U_0}{\partial x}=-\lambda_{01}U_0+\sum_{i=1}^s\lambda_{i1}U_i+\mathcal{O}(h_i^2),\\[2mm]
		&\frac{\partial^2 U_0}{\partial x^2}=-\lambda_{02}U_0+\sum_{i=1}^s\lambda_{i2}U_i+\mathcal{O}(h_i^2),\\[2mm]
	\end{aligned}
	\right.
\end{equation}
\section{2D explicit formulae}
Let $\Omega=[0,L]\times[0,L] \subset \mathbb{R}^2$ be a domain and $$M= \{\textbf{x}_1, \dots, \textbf{x}_N \} \subset \Omega$$ a discretization of $\Omega$ with $N$ points (see Figure 1). 
For each one of the nodes of the domain, where the value of $u$ is unknown, a star is defined as a set of selected points $$\{\textbf{x}_i:i=1,...,s,1\leq i\leq s\} \subset M$$ with the central node $\textbf{x}_0 \in M$ and $\textbf{x}_i$, $(i=1,\dots,s) \in M$ is a set of points located in the neighbourhood of $\textbf{x}_0$. In order to select the points different criteria as four quadrants or distance can be used \cite{R1}.
\\
Let $\textbf{x}_0=(x_0,y_0)$ be the central node of a star and $h_i=x_i-x_0, k_i=y_i-y_0$, where $(x_i,y_i)$ are the coordinates of the $i^{th}$ node of the star. Then by the Taylor series expansion we have
\begin{equation}
	u_i = u_0 + h_i \dfrac{\partial u_0}{\partial x} + k_i \dfrac{\partial u_0}{\partial y}+ \dfrac{1}{2} \bigg( h^2_i \dfrac{\partial^2 u_0}{\partial x^2} + k^2_i \dfrac{\partial^2 u_0}{\partial y^2} + 2 h_i k_i \dfrac{\partial^2 u_0}{\partial x \partial y} \bigg) + ...,  \label{explicito1}\end{equation}
for $i=1,...,s$.
\\
Let us use the notations
\begin{equation*}{\textbf{c}_i}^T=\{h_i,k_i,\frac{h_i^2}{2},\frac{k_i^2}{2},h_ik_i\}
\end{equation*}
and
\begin{equation*}\boldsymbol{D_5}^T=\{\frac{\partial u_0}{\partial x},\frac{\partial u_0}{\partial y},\frac{\partial^2u_0}{\partial x^2},\frac{\partial^2u_0}{\partial y^2},\frac{\partial^2u_0}{\partial x \partial y}\}.\end{equation*}
If  we do not consider in (\ref{explicito1}) the higher than second order terms, we can obtain a second order approximation of $u_i$, which we shall denote $U_i$. Then, we define  the following:
\begin{equation}
	\begin{array}{ll}
		B(U)&\displaystyle=\displaystyle\sum_{i=1}^s[(U_0-U_i)+h_i\frac{\partial U_0}{\partial x}+k_i\frac{\partial U_0}{\partial y}+
		\\\\
		&\displaystyle+\frac{1}{2}(h_i^2\frac{\partial^2U_0}{\partial x^2}+k_i^2\frac{\partial^2U_0}{\partial y^2}+2h_ik_i\frac{\partial^2U_0}{\partial x \partial y})]^2w_i^2\label{explicito2},
	\end{array}
\end{equation}
where $w_i=w(h_i,k_i)$ are positive symmetrical weighting functions decreasing in magnitude as the distance to the center increases, as defined in Lankaster and Salkauskas \cite{lancaster} (see also Levin \cite{R20}). Another weighting functions as potentials or exponential can be used (see \cite{R6} for more details). We minimize the norm given by (\ref{explicito2}) with respect to the partial derivatives by considering the following linear system
\begin{equation*} \boldsymbol{A}(h_i,k_i,w_i)\boldsymbol{D_5}=\boldsymbol{ b}(h_i,k_i,w_i,U_0,U_i)\end{equation*}
where
\begin{equation*}
	\boldsymbol{A}=\left(
	\begin{array}{cccc}
		h_1 & h_2 & \cdots & h_s \\
		k_1 & k_2 & \cdots & k_s \\
		\vdots & \vdots & \vdots & \vdots \\
		h_1k_1 & h_2k_2 & \cdots & h_sk_s \\
	\end{array}
	\right)\left(
	\begin{array}{cccc}
		\omega_1^2 &  &  &  \\
		& \omega_2^2 &  &  \\
		&  & \cdots &  \\
		&  &  & \omega_s^2 \\
	\end{array}
	\right)\left(
	\begin{array}{cccc}
		h_1 & k_1 & \cdots & h_1k_1 \\
		h_2 & k_2 & \cdots & h_2k_2 \\
		\vdots & \vdots& \vdots & \vdots \\
		h_s & k_s &  \cdots & h_sk_s \\
	\end{array}\right),\end{equation*}
and
\begin{equation*}
	\begin{array}{ll}
		\boldsymbol{ b}^T=&\displaystyle \left(\sum_{i=1}^{s}(-U_0+U_i)h_iw_i^2,\;\sum_{i=1}^{s}(-U_0+U_i)k_iw_i^2,\;\sum_{i=1}^{s}(-U_0+U_i)\frac{h_i^2w_i^2}{2},\right.\\
		&\displaystyle
		\left.\sum_{i=1}^{s}(-U_0+U_i)\frac{k_i^2w_i^2}{2},\;\sum_{i=1}^{s}(-U_0+U_i)h_ik_iw_i^2\right).
	\end{array}
\end{equation*}
It is well known that  $\boldsymbol{A}$ is a positive definite matrix and the approximation is of second order $\Theta(h_i^2,k_i^2)$ (see \cite{R6}). If we define
\begin{equation*} \boldsymbol{A}^{-1}=\boldsymbol{Q}\boldsymbol{Q}^T,\end{equation*}
we have
\begin{equation} \boldsymbol{D_5}=\boldsymbol{Q}\boldsymbol{Q}^T\boldsymbol{b}.\label{explicito3}\end{equation}
Thus, equation (\ref{explicito3}) can be rewritten as
\begin{equation*} \boldsymbol{D_5}=-U_0\boldsymbol{Q}\boldsymbol{Q}^T\displaystyle\sum_{i=1}^s w_i^2\boldsymbol{c_i}+\boldsymbol{Q}\boldsymbol{Q}^T\displaystyle\sum_{i=1}^sU_i w_i^2\boldsymbol{c_i},\end{equation*}
or
\begin{equation*} \boldsymbol{D}=\boldsymbol{Q}\boldsymbol{Q}^T\boldsymbol{W}(\boldsymbol{u}-u_0\boldsymbol{1})\end{equation*}
where
$$\boldsymbol{W}=\left(
\begin{array}{cccc}
	h_1w_1^2 & h_2w_2^2 & \cdots & h_sw_s^2 \\
	k_1w_1^2 & k_2w_2^2 & \cdots & k_sw_s^2 \\
	\frac{h_1^2}{2}w_1^2 &  \frac{h_2^2}{2}w_2^2 & \vdots &  \frac{h_s^2}{2}w_s^2\\
	\frac{k_1^2}{2}w_1^2 &  \frac{k_2^2}{2}w_2^2 & \vdots &  \frac{k_s^2}{2}w_s^2\\
	h_1k_1w_1^2 & h_2k_2 w_s^2& \cdots & h_sk_sw_s^2 \\
\end{array}
\right)$$
and
$$\boldsymbol{1}=\left\{1,1,\cdots,1\right\};\quad \boldsymbol{U}=\left\{U_1,U_2,\cdots,U_s\right\}^T.$$
As in \cite{R7}, we denote the spatial derivatives using GFD   by
\begin{equation}\left\{\begin{array}{l}
		\dfrac{\partial U(\textbf{x}_0,n\Delta t)}{\partial x}= -m_{01}U_0+\displaystyle\sum_{i=1}^s m_{i1}U_i+\Theta(h_i^2,k_i^2),\: with\: m_{01}=\displaystyle\sum_{i=1}^s m_{i1}
		\\\dfrac{\partial U(\textbf{x}_0,n\Delta t)}{\partial y}= -m_{02}U_0+\displaystyle\sum_{i=1}^s m_{i2}U_i+\Theta(h_i^2,k_i^2),\: with\: m_{02}=\displaystyle\sum_{i=1}^s m_{i2}
		\\
		\Biggl(\dfrac{\partial^2 U}{\partial x^2}+\dfrac{\partial^2 U}{\partial y^2}\Biggr)|_{(\textbf{x}_0,n\Delta t)}= -m_{00}U_0+\displaystyle\sum_{i=1}^s m_{i0}U_i+\Theta(h_i^2,k_i^2),\\ m_{00}=\displaystyle\sum_{i=1}^s m_{i0},\\
	\end{array}\right.\label{explicito4}
\end{equation}
or, written in vectorial form,
\begin{equation*}
	\boldsymbol{D_5}U(\textbf{x}_0,n\Delta t)=-\boldsymbol{m_0}U_0+\sum^{s}_{i=1}\boldsymbol{m_i}U_i+\Theta(h_i^2,k_i^2)\end{equation*}
where $\boldsymbol{m_0}$ and $\boldsymbol{m_i}$ stand for
\begin{equation*}\begin{split}
		\boldsymbol{m_0}=\{m_{01},m_{02},m_{03},m_{04},m_{05}\}^{T},\\
		\boldsymbol{m_i}=\{m_{i1},m_{i2},m_{i3},m_{i4},m_{i5}\}^{T},\\
		m_{00}=m_{03}+m_{04};\quad m_{i0}=m_{i3}+m_{i4},\\
	\end{split}
\end{equation*}
fulfilling
\begin{equation*}
	\boldsymbol{m_0}=\sum_{i=1}^{s}\boldsymbol{m_i}.
\end{equation*}
\section{GFDM scheme}\label{s3}
We present the explicit GFDM scheme for the one and two dimensional case. Consider first $\Omega\subset\mathbb{R}$. We use the formulae (\ref{discretizacion}) for the approximation of the spatial derivatives and the time derivative approximation
\begin{equation}\label{timeApprox}
	\frac{\partial U(x_0,n\Delta t)}{\partial t}=\frac{U^{n+1}_0-U^{n}_0}{\Delta t}+\mathcal{O}(\Delta t),\quad \Omega\subset\mathbb{R}
\end{equation}
Hence, our numerical scheme is
\begin{equation}
	\label{esquema1d}
	\left\lbrace
	\begin{split}
		&\dfrac{k^{n+1}_0-k^{n}_0}{\Delta t}=-\lambda_{00}k^n_0+\sum_{i=1}^{s}\lambda_{i0}k^n_i+^{1D}\mathbb{F}^n_{0,i}+A^n_0f(k^n_0)-\delta k^n_0+\mathcal{O}(\Delta t,h_i^2)\\
		&\dfrac{A^{n+1}_0-A^{n}_0}{\Delta t}=d\Biggl(-\lambda_{00}A^n_0+\sum_{i=1}^{s}\lambda_{i0}A^n_i\Biggr)+A^n_0g(A^n_0)+\mathcal{O}(\Delta t,h_i^2),
	\end{split}
	\right.
\end{equation}
where 
\begin{equation*}
	\begin{split}
		^{1D}\mathbb{F}^n_{0,i}:=&-\chi\Biggl(-\lambda_{01}k^n_0+\sum_{i=1}^s\lambda_{i1}k^n_i\Biggr)\Biggl(-\lambda_{01}A^n_0+\sum_{i=1}^s\lambda_{i1}A^n_i\Biggr)\\&-\chi k^n_0\Biggl(-\lambda_{00}A^n_0+\sum_{i=1}^{s}\lambda_{i0}A^n_i\Biggr).
	\end{split}
\end{equation*}
For the 2-dimensional case, the GFD explicit scheme is
\begin{equation}
	\label{esquema2d}
	\left\lbrace
	\begin{split}
		&\dfrac{k^{n+1}_0-k^{n}_0}{\Delta t}=-m_{00}k^n_0+\sum_{i=1}^{s}m_{i0}k^n_i+^{2D}\mathbb{F}^n_{0,i}+A^n_0f(k^n_0)-\delta k^n_0+\mathcal{O}(\Delta t,h_i^2,k_i^2)\\
		&\dfrac{A^{n+1}_0-A^{n}_0}{\Delta t}=d\Biggl(-m_{00}A^n_0+\sum_{i=1}^{s}m_{i0}A^n_i\Biggr)+A^n_0g(A^n_0)+\mathcal{O}(\Delta t,h_i^2,k_i^2),
	\end{split}
	\right.
\end{equation}
where \begin{equation*}
	\begin{split}
		^{2D}\mathbb{F}^n_{0,i}:=&-\chi\Biggl[\Biggl(-m_{01}k^n_0+\sum_{i=1}^sm_{i1}k^n_i\Biggr)\Biggl(-m_{01}A^n_0+\sum_{i=1}^sm_{i1}A^n_i\Biggr)+\\
		&+\Biggl(-m_{02}k^n_0+\sum_{i=1}^sm_{i2}k^n_i\Biggr)\Biggl(-m_{02}A^n_0+\sum_{i=1}^sm_{i2}A^n_i\Biggr)\Biggr]\\&-\chi k^n_0\Biggl(-m_{00}A^n_0+\sum_{i=1}^{s}m_{i0}A^n_i\Biggr).
	\end{split}
\end{equation*}
\begin{rem}
	The consistency of the GFDM formulae was proved in \cite{vargas} in 1D and in \cite{R6} in 2D.
\end{rem}
We provide the proof of the convergence of the GFD scheme to the continuous solution of the PDE system for the two dimensional case only since the one dimensional case is identical. The result states as follows:
\begin{thm}\label{thm}
	Let $k,A$ be the exact solution of (\ref{1}). Then, the GFD explicit scheme (\ref{esquema2d}) is convergent if
	\begin{equation}\label{hip1}
		0<m_{00}+\Phi_1-\Phi_2
	\end{equation}
	and 
	\begin{equation}\label{hip2}
		\Delta t<\dfrac{2}{m_{00}+\Phi_1+\Phi_2}
	\end{equation}
	for $\Phi_1$ and $\Phi_2$ explicitly given in the proof.
\end{thm}
{\bf Proof of Theorem \ref{thm}}
\newline
We take the difference between GFD scheme (\ref{esquema2d}) and the expression for the exact solution. Let $\overline{k}_j^n$ be the approximated $k$-solution at time $n\Delta t$ (similarly $\overline{A}^n_j$) and $k^n_j$ the value of the exact  $k$-solution (similarly $A^n_j$). Also, we call $\tilde{k}_j^n=k_j^n-\overline{k}_j^n\textnormal{, }\tilde{A}_j^n=A_j^n-\overline{A}_j^n$. 
$$\dfrac{A^{n+1}_0-A^{n}_0}{\Delta t}=A^{n}_0g(A^n_0)+\mathcal{O}(\Delta t,h_i^2,k_i^2).$$
Let us define $\tilde{k}^n=\displaystyle\max_{i\in\{0,...,s\}}|\tilde{k}_i^n|$ and $\tilde{A}^n=\displaystyle\max_{i\in\{0,...,s\}}|\tilde{A}_i^n|$. Therefore the equation for the error $\tilde{A}^n$ reads
$$\tilde{A}^{n+1}\leq \tilde{A}^{n}\cdot\rho+\tau_i,$$
where $\rho=|1+\Delta t(g(A^n)+A^ng'(\eta))|$ and $\tau_i=\mathcal{O}(\Delta t(\Delta t,h_i^2,k_i^2))$
By induction, we see that
$$\tilde{A}^1=\tilde{A}^0\rho+\tau_i=\tau_i,$$
$$\tilde{A}^2=\tilde{A}^1\rho+\tau_i=\tau_i(1+\rho)$$
and 
$$\tilde{A}^{n+1}=\tau_i(1+\rho+...+\rho^n).$$
Now, let us call $T$ the final time where we obtain the numerical solution so $n\Delta t\leq T$. Then, 
$$\tilde{A}^{n+1}\leq \tau_i\sum_{k=0}^n\rho^k\leq \tau_ie^{T|(g(A^n)+A^ng'(\eta)|}.$$
Finally, $\tilde{A}^{n+1}=\mathcal{O}(\Delta t(\Delta t,h_i^2,k_i^2))$ as $t,h_i,k_i\to 0$.
Now, the following expression yields
\begin{equation}
	\begin{split}\label{esquema2d7}
		\tilde{k}^{n+1}_0=&\tilde{k}^{n}_0\Biggl[1-\Delta t(m_{00}+\delta)\Biggr]+\Delta t [A^{n}_0f(k^{n}_0)-\overline{A}^{n}_0f(\overline{k}^{n}_0)]\\&+\Delta t\sum_{i=1}m_{i0}\tilde{k}^n_i+\Delta t(\mathbb{F}^n_{0,i}-\mathbb{\overline{F}}^n_{0,i})+\mathcal{O}(\Delta t(\Delta t,h_i^2,k_i^2)).
	\end{split}
\end{equation}
Now, we use the mean value theorem:
\begin{equation}\label{esquema2d8}
	\begin{split}
		A^{n}_0f(k^{n}_0)-\overline{A}^{n}_0f(\overline{k}^{n}_0)=&A^{n}_0f(k^{n}_0)\pm\overline{A}^{n}_0f(k^{n}_0)-\overline{A}^{n}_0f(\overline{k}^{n}_0)\\
		&=\tilde{A}^{n}_0f(k^{n}_0)+\overline{A}^{n}_0f'(\xi)\tilde{k}^{n}_0,
	\end{split}
\end{equation}
for some $\xi\in(k^{n}_0,\overline{k}^{n}_0)\cup(\overline{k}^{n}_0,k^{n}_0)$. Then, substituting in \ref{esquema2d7}
\begin{equation}
	\begin{split}\label{esquema2d10}
		\tilde{k}^{n+1}_0=&\tilde{k}^{n}_0\Biggl[1-\Delta t(m_{00}+\delta-\overline{A}^{n}_0f'(\xi)\tilde{k}^{n}_0)\Biggr]+\Delta t \tilde{A}^{n}_0f(k^{n}_0) \\&+\Delta t\sum_{i=1}^s m_{i0}\tilde{k}^n_i+\Delta t(\mathbb{F}^n_{0,i}-\mathbb{\overline{F}}^n_{0,i})+\mathcal{O}(\Delta t(\Delta t,h_i^2,k_i^2)).
	\end{split}
\end{equation}
Lets look at the term $\Delta t(\mathbb{F}^n_{0,i}-\mathbb{\overline{F}}^n_{0,i})$. First, 
$$-\chi(m^2_{01}k^n_0A^n_0-m^2_{01}\overline{k}^n_0\overline{A}^n_0\pm m^2_{01}k^n_0\overline{A}^n_0)=-\chi(m^2_{01}k^n_0\tilde{A}^n_0+m^2_{01}\tilde{k}^n_0\overline{A}^n_0).$$
Second,
\begin{equation*}
	\begin{split}
		\chi\Biggl[m_{01}k^n_0&\sum_{i=1}^sm_{i1}A^n_i-m_{01}\overline{k}^n_0\sum_{i=1}^sm_{i1}\overline{A}^n_i\pm m_{01}k^n_0\sum_{i=1}^sm_{i1}\overline{A}^n_i\Biggr]\\&=\chi\Biggl[m_{01}k^n_0\sum_{i=1}^sm_{i1}\tilde{A}^n_i+m_{01}\tilde{k}^n_0\sum_{i=1}^sm_{i1}A^n_i\Biggr]
	\end{split}
\end{equation*}
For the third term,
\begin{equation*}
	\begin{split}
		\chi\Biggl[m_{01}A^n_0&\sum_{i=1}^sm_{i1}k^n_i-m_{01}\overline{A}^n_0\sum_{i=1}^sm_{i1}\overline{k}^n_i\pm m_{01}\overline{A}^n_0\sum_{i=1}^sm_{i1}k^n_i\Biggr]\\&=\chi\Biggl[m_{01}\tilde{A}^n_0\sum_{i=1}^sm_{i1}k^n_i+m_{01}\overline{A}^n_0\sum_{i=1}^sm_{i1}\tilde{k}^n_i\Biggr]
	\end{split}
\end{equation*}
For the last term,
\begin{equation*}
	\begin{split}
		-\chi\Biggl[\sum_{i=1}^sm_{i1}k^n_i&\sum_{i=1}^sm_{i1}A^n_i-\sum_{i=1}^sm_{i1}\overline{k}^n_i\sum_{i=1}^sm_{i1}\overline{A}^n_i\pm \sum_{i=1}^sm_{i1}k^n_i\sum_{i=1}^sm_{i1}\overline{A}^n_i\Biggr]\\&=\chi\Biggl[\sum_{i=1}^sm_{i1}k^n_i\sum_{i=1}^sm_{i1}\tilde{A}^n_i+\sum_{i=1}^sm_{i1}\tilde{k}^n_i\sum_{i=1}^sm_{i1}\overline{A}^n_i\Biggr]
	\end{split}
\end{equation*}
Hence, putting all together,
\begin{equation}
	\begin{split}\label{esquema2d11}
		\tilde{k}^{n+1}\leq&\tilde{k}^{n}\Biggl[\Biggl|1-\Delta t\Biggl(m_{00}+\delta-\overline{A}_0^{n}f'(\xi)+\chi m^2_{01}\overline{A}^n_0+\chi m^2_{02}\overline{A}^n_0\\&+\chi m_{01}\sum_{i=1}^sm_{i1}A^n_i+\chi m_{02}\sum_{i=1}^sm_{i2}A^n_i\\&-\chi\Biggl(-m_{00}\overline{A}^n_0+\sum_{i=1}^sm_{i0}\overline{A}^n_i\Biggr)\Biggr)\Biggr|+\Delta t\Biggl(\sum_{i=1}^s|m_{i0}|\\&
		+|\chi m_{01}\overline{A}^n_0|\sum_{i=1}^s|m_{i1}|+|\chi m_{02}\overline{A}^n_0|\sum_{i=1}^s|m_{i2}|\\&+|\chi|\sum_{i=1}^s|m_{i1}||\sum_{i=1}^sm_{i1}\overline{A}^n_i|+|\chi|\sum_{i=1}^s|m_{i2}||\sum_{i=1}^sm_{i2}\overline{A}^n_i|\Biggr)\Biggr]\\&+\Delta t \tilde{A}^{n}\Biggl[\Biggl|f(k^{n})-\chi m^2_{01}k^n_0-\chi m^2_{02}k^n_0\\&+\chi m_{01}\sum_{i=1}^sm_{i1}k^n_i+\chi m_{02}\sum_{i=1}^sm_{i2}k^n_i-\chi m_{00}k^n_0\Biggr|\\&+|\chi m_{01}k^n_0|\sum_{i=1}^s|m_{i1}|+|\chi m_{02}k^n_0|\sum_{i=1}^s|m_{i2}|\\&+|\chi\sum_{i=1}^sm_{i1}k^n_0|\sum_{i=1}^s|m_{i1}|+|\chi\sum_{i=1}^sm_{i2}k^n_0|\sum_{i=1}^s|m_{i2}|\\&+|\chi k^n_0|\sum_{i=1}^s|m_{i0}|\Biggr] +\mathcal{O}(\Delta t(\Delta t,h_i^2,k_i^2)).
	\end{split}
\end{equation}
We can write the last inequality as 
\begin{equation}
	\tilde{k}^{n+1}\leq\tilde{k}^{n}|1-\Delta t (m_{00}+\Phi_1)|+\Delta t\Phi_2,
\end{equation}
for an obvious choice of $\Phi_1$ and $\Phi_2$. Convergence follows from
\begin{equation}
	|1-\Delta t (m_{00}+\Phi_1)|+\Delta t\Phi_2<1.
\end{equation}
First, 
$$-1+\Delta t\Phi_2<1-\Delta t(m_{00}+\Phi_1),$$
true by (\ref{hip1})
and second,
$$1-\Delta t (m_{00}+\Phi_1)<1-\Delta t\Phi_2,$$
which holds by (\ref{hip2}).
$\Box$
\begin{rem}
	The inequalities  of (\ref{hip1}) and (\ref{hip2}) give us a range of values for $\Delta t$ for convergence of each one of the stars of the domain. Then the minimum value obtained among all the stars is taken as $\Delta t$ for convergence condition.
\end{rem}
\section{Numerical examples}
In this section we present several 1D and 2D examples where the dynamics of the numerical model are shown. We use the dimensionless form of the model of \cite{capasso}. For the 1D examples we use stars with 2 nodes whereas for the 2D examples 8-node stars are employed. In both cases, we use $\Delta t=0.001$.
\subsection{1D examples}
\begin{figure}[h]
	\begin{center}
		\quad \includegraphics[width=0.45\textwidth]{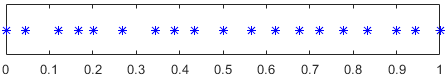}
		\caption{1D irregular clouds of points }
		\label{figure1}
	\end{center}
\end{figure}

First, we look at the one dimensional case and use the irregular discretization of $[0,1]$ of Figure \ref{figure1}. 
The initial data for $k$ is always chosen as
$$k_0(x)=\left\lbrace
\begin{aligned}{}
	&5& x\in [0,0.25], \\
	&40x-5 & \quad x\in (0.25,0.75), \\
	&25 &  x\in [0.75,1].  \\
\end{aligned}
\right.$$
and 
$$A_0(x)=1,\quad g(x)=0.1e^{\frac{-(x-0.5)^2}{2\cdot 0.2^2}}.$$
\\

\begin{figure}[h]
	\begin{center}
		\includegraphics[width=0.40\textwidth]{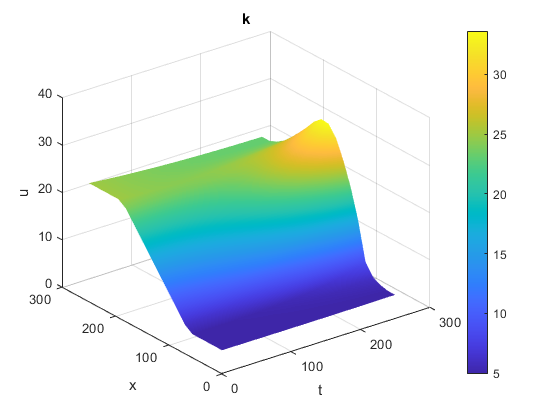}
		\includegraphics[width=0.40\textwidth]{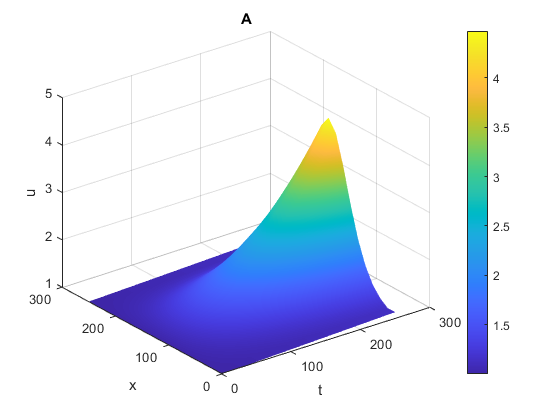}
		\caption{Solution $k$ for $\delta=0.05$.}\label{figure13}
	\end{center}
\end{figure}

\begin{figure}[h]
	\begin{center}
		\includegraphics[width=0.45\textwidth]{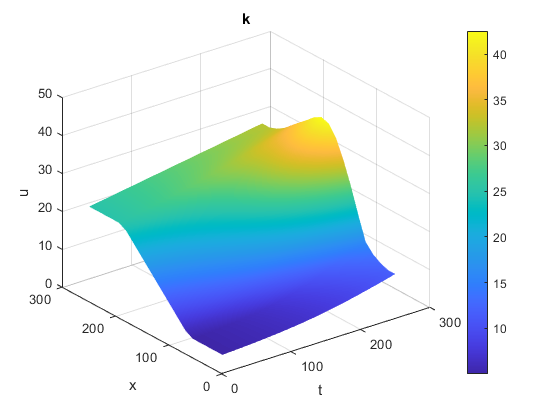}
		\includegraphics[width=0.45\textwidth]{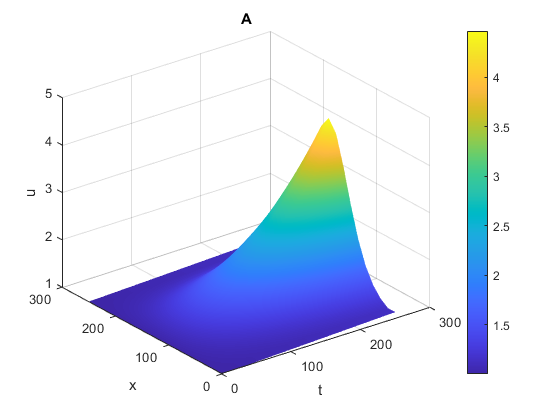}
		\caption{Solution $k$ for $\delta=0.02$.}\label{figure14}
	\end{center}
\end{figure}

\begin{figure}[h]
	\begin{center}
		\includegraphics[width=0.45\textwidth]{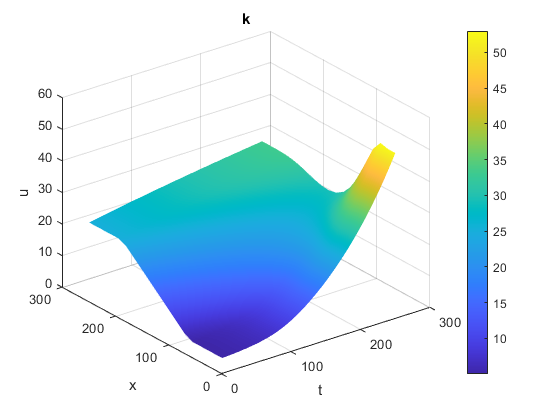}
		\includegraphics[width=0.45\textwidth]{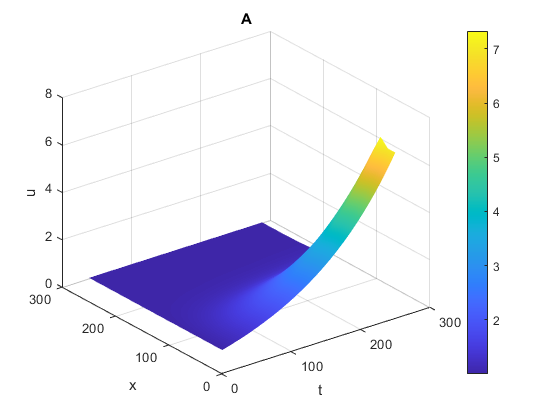}
		\caption{Solution $k$ for $\delta=0.02$.}\label{figure15}
	\end{center}
\end{figure}

We start investigating the influence of depreciation in the model, and choose:
Figure \ref{figure13}
$\delta=0.05$, $L=20$, $A_0(x)=1$, $g(x)=0.1e^{\frac{-(x-0.5)^2}{2\cdot 0.2^2}}$. 
\\
By choosing now a lower depreciation rate, $\delta=0.02$, and the same relation of parameters, the numerical solution shows how the capital distribution grows at every point of the domain (Figure \ref{figure14}). 
\\
We consider, next, the effects of the capital movement towards the regiosn with high concentration of technology. To further explore this idea, we consider a different regional technological growth rate and put
$$ g(x)=0.1e^{\frac{-(x-0.1)^2}{2\cdot 0.2^2}},\quad \chi=1.$$
In this way, we model the situation in which the technological production comes from a poor region. The dynamics of the system shows a fast capital growth, particularly in the poorest regions. The result can be seen in Figure \ref{figure15}. Notice that the induced taxis term provokes an reinforment in the distribution in the sense that the maximum value for $k$ (50, approximately) is greatest than in the previous case (40).
\clearpage
\subsection{2D examples}
\begin{figure}[h]
	\begin{center}
		\includegraphics[width=0.50\textwidth]{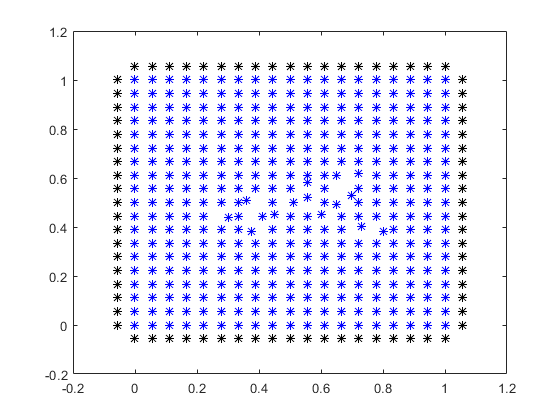}
		\caption{2D irregular cloud of points.}\label{figura2d}
	\end{center}
\end{figure}

For the 2D case, we employ the irregular discretization of $[0,1]\times[0,1]$ of Figure \ref{figura2d}.	
In all examples, initial data is chosen as in Figure 6.
\begin{figure}[h]
	\begin{center}
		\includegraphics[width=0.50\textwidth]{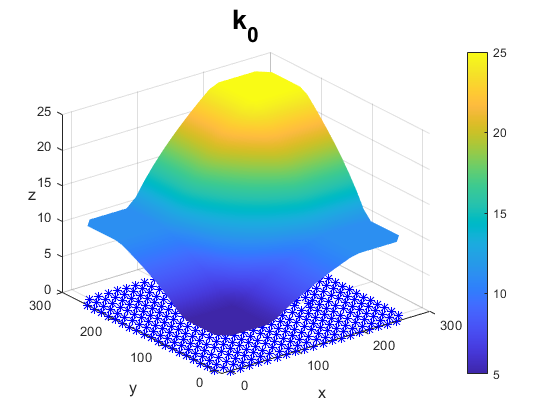}
		\label{dato2d}
		\caption{Initial data for $k_0$.}
	\end{center}
\end{figure}
\\
\\
First, we investigate the influence of the depreciation rate in the two-dimensional model with non-constant technological progress. We set $\mu=0.05$, $L=150$, and a technology growth rate of $g(x)=0.1e^{\frac{-(x-0.5)^2-(y-0.5)^2}{2\cdot 0.2^2}}$, the initial technology distribution $A_0(x)=1$ and $\chi=0$ (Figure \ref{figure16}). 
We plot the numerical solution and observe that the per capita capital distribution takes on the shape of the technological progress.
\begin{figure}[h]
	\begin{center}
		\includegraphics[width=0.45\textwidth]{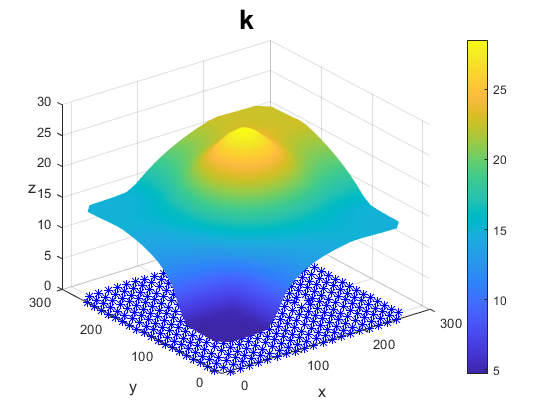}
		\includegraphics[width=0.45\textwidth]{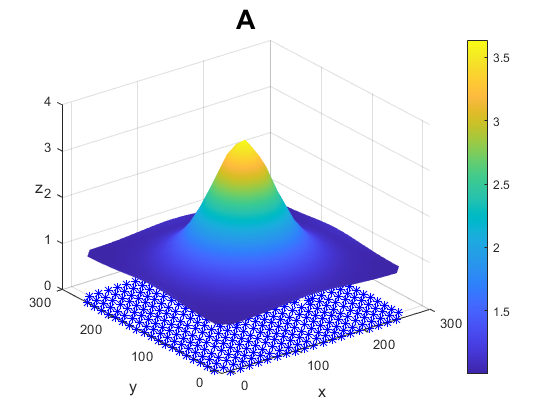}
		\caption{Solution $k$ for $\delta=0.05$.}\label{figure16}
	\end{center}
\end{figure}
\\
In Figure \ref{figure17}, we model the poverty trap in the 2D case for a constant $A=1$ and $\mu=0.085$. At an intermediate time, the richest regions seem to approach a positive steady state, but the high depreciation rate causes them to ultimately converge to zero. By assuming non-constant technological progress and using the same range of parameters, Figure \ref{figure18} shows that although the capital density initially decreases, it eventually approaches the shape of technology.
\begin{figure}[h]
	\begin{center}
		\includegraphics[width=0.45\textwidth]{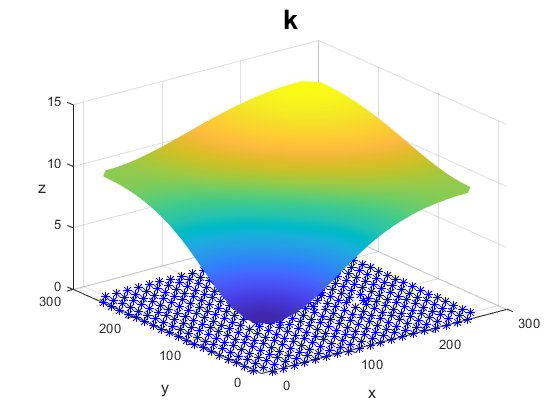}
		\includegraphics[width=0.45\textwidth]{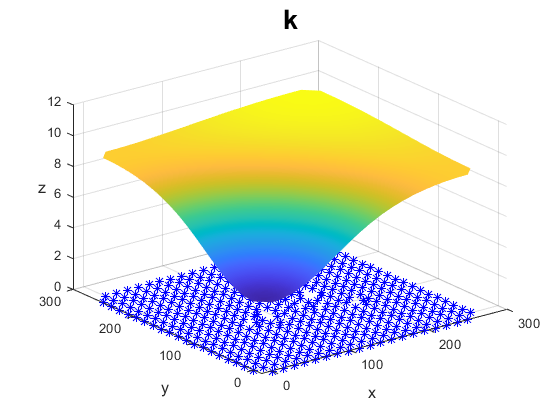}\\
		\includegraphics[width=0.45\textwidth]{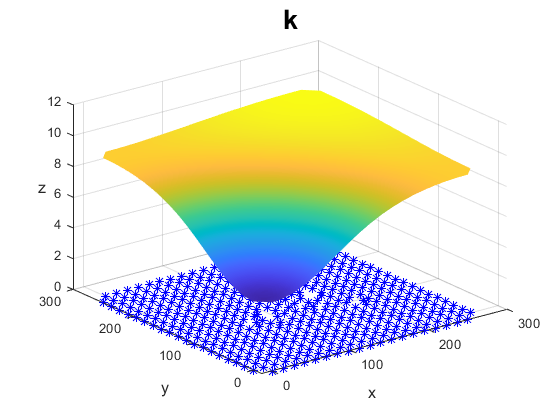}
		\includegraphics[width=0.45\textwidth]{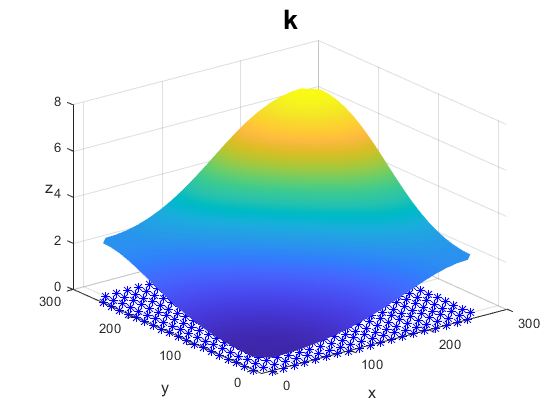}
		\caption{Solution $k$ for $\delta=0.085$.}\label{figure17}
	\end{center}
\end{figure}
\begin{figure}[h]
	\begin{center}
		\includegraphics[width=0.45\textwidth]{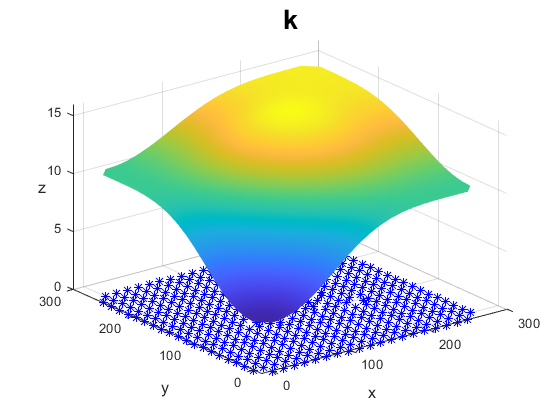}
		\includegraphics[width=0.45\textwidth]{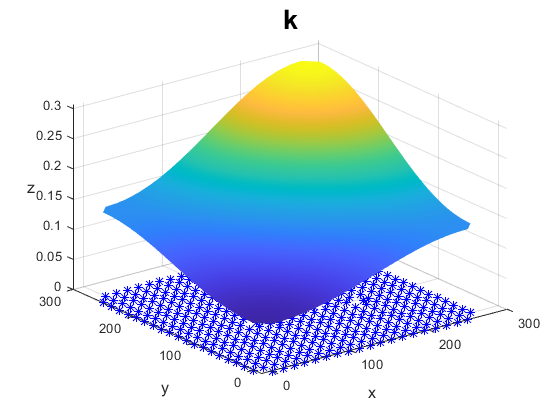}\\
		\includegraphics[width=0.45\textwidth]{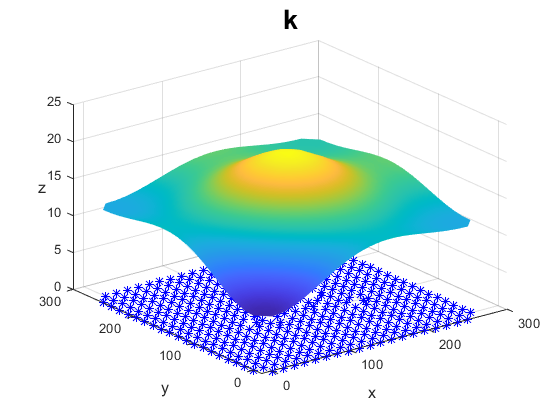}
		\includegraphics[width=0.45\textwidth]{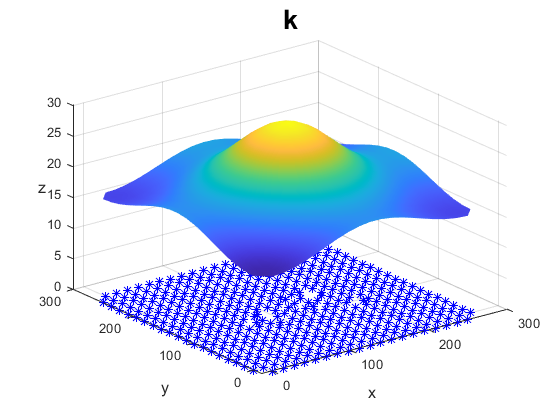}
		\caption{Solution $k$ for $\delta=0.085$.}\label{figure18}
	\end{center}
\end{figure}
\\
\\
Lastly, we showcase an illustration featuring technology-induced mobility and an exceptionally high depreciation rate of $\mu=0.3$. When $\chi=0$, the convergence towards zero happens immediately. However, when we consider the movement of capital towards regions with the most advanced technology, the capital density eventually rebounds over time, potentially producing solution blow-ups. It is worth noting that, unlike the other cases, the per capita capital exhibits a spiky distribution, which models a scenario where inequalities tend to escalate during a period of high depreciation rate (see Figure \ref{figure19}).
\begin{figure}[h]
	\begin{center}
		\includegraphics[width=0.45\textwidth]{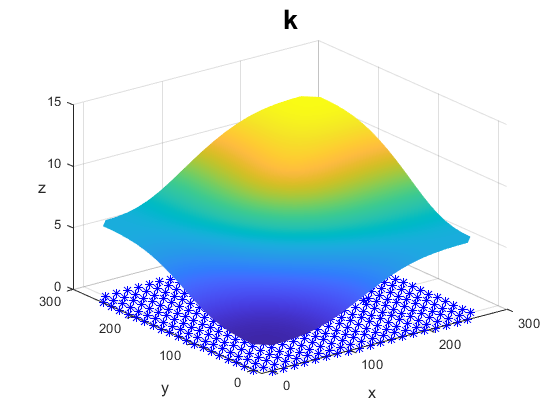}
		\includegraphics[width=0.45\textwidth]{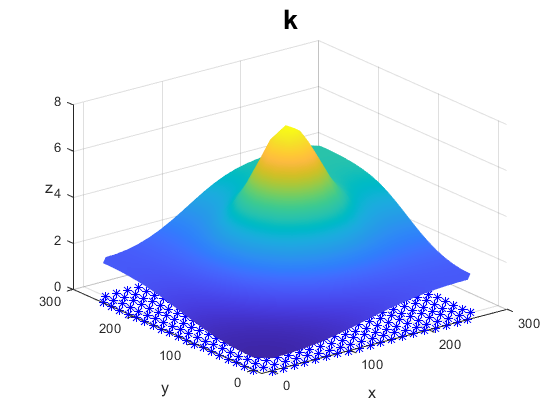}\\
		\includegraphics[width=0.45\textwidth]{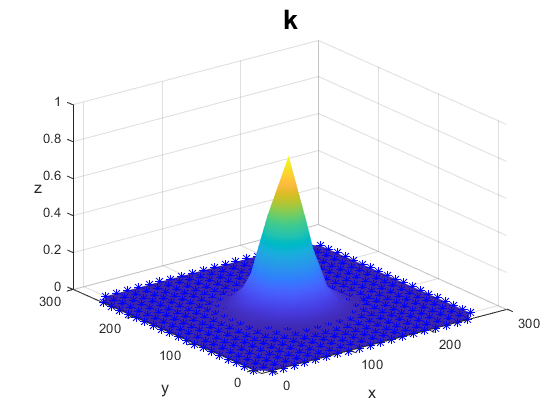}
		\includegraphics[width=0.45\textwidth]{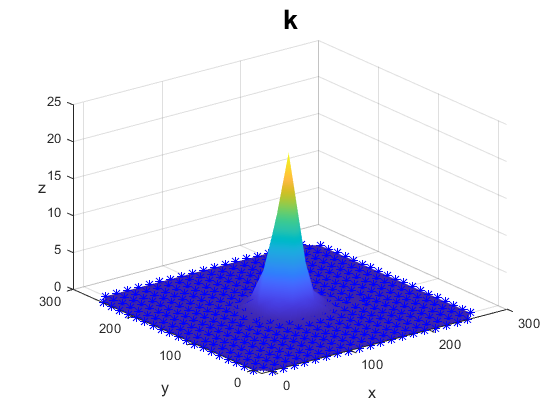}
		\caption{Solution $k$ for $\delta=0.3$.}\label{figure19}
	\end{center}
\end{figure}
\section{Conclusions}
We have employed the Generalized Finite Difference Method (GFDM) to derive the discretization of a system of Partial Differential Equations (PDEs) described by equation (\ref{1}). In Theorem \ref{thm}, we have established the conditional convergence of this method for solving the nonlinear system and explicitly provided the convergence condition.
\\
However, the model does have a notable limitation. The equation governing technological progress remains unaffected by capital, potentially leading to the complete disappearance of capital without impeding technological advancement. Additionally, the model assumes a constant available workforce. To overcome these challenges and gain a deeper understanding, further investigation using analytical and numerical techniques is essential.
\\
To validate the asymptotic behavior of the solution stated in the theory and to demonstrate the accuracy and efficiency of GFDM applied to this highly nonlinear system of coupled parabolic PDEs over irregular domains, numerical tests are proposed.

\section*{Acknowledgements}
AMV is supported by the Spanish MINECO through Juan de la Cierva fellow-ship FJC2021-046953-I.

%
%



%
%
%
\end{document}